# Making proofs without Modus Ponens: An introduction to the combinatorics and complexity of cut elimination

A. Carbone and S. Semmes [1]

**Abstract**

Modus Ponens says that if you know $A$ and you know that $A$ implies $B$ then you know $B$. This is a basic rule that we take for granted and use repeatedly, but there is a gem of a theorem in logic by Gentzen to the effect that it is not needed in some logical systems. It is fun to say "You can make proofs without lemmas" to mathematicians and watch how they react, but our true intention here is to let go of logic as a reflection of reasoning and move towards combinatorial aspects. Proofs contain basic problems of algorithmic complexity within their framework, and there is strong geometric and dynamical flavor inside them.

## 1 The beginning of the story

Mathematicians are making proofs every day. In proof theory one studies proofs. This is frightening to many mathematicians, but a principal theme of the present exposition is to treat logic unemotionally.

The idea of complexity sheds an interesting light on proofs. A basic question is whether a propositional tautology of size $n$ should always have a short proof, a proof of size $p(n)$ for some polynomial $p$ for instance. There is a proof system in which this is true if and only if "$NP = co - NP$". The latter is an unsolved general question in computational complexity which is related to the existence of polynomial-time algorithms in cases where only exponential algorithms are known. We shall say more about this later. The equivalence was established in [16].

Sometimes proofs have to be long if one does not permit a rule like Modus Ponens. Such a rule allows *dynamics* within the implicit computations occurring in proofs. In tracing through a proof one may visit the same formula repeatedly with many substitutions. The level in the hierarchy of proof systems at which dynamics appears is the first at which little is known concerning the existence of short proofs of hard tautologies. At lower levels it is known that there are tautologies with only exponential size proofs. Dynamical structure in proofs seems to be important for making proofs short. (This is discussed further in [14].)

---

[1] The first author was supported by the Lise-Meitner Stipendium # M00187-MAT (Austrian FWF) and the second author was supported by the U.S. National Science Foundation. Both authors are grateful to IHES for its hospitality.

The idea of dynamical structure in proofs may seem odd at first, but it is very natural, both in ordinary mathematical activity and in formal logic. There is a way to eliminate the nontrivial dynamical structure from a proof, but at great cost in expansion. This is the matter of cut elimination that we shall discuss here.

Before we get to that we should review some background information about logic and complexity. Our discussion will be informal, but one can find more precision and detail in [30, 36, 48]. We begin with some comments about propositional logic, predicate logic, and arithmetic.

Propositional logic is the simplest of the three. One has variables, often called $p$, $q$, etc., and one can build formulas out of them using standard connectives, $\vee$ (or), $\wedge$ (and), $\neg$ (negation). In predicate logic one uses the same connectives, and one also has the quantifiers $\forall$ (for all) and $\exists$ (there exists) for making formulas like $\exists x F(x)$. Here $x$ is a variable and $F$ is a relation, in this case a unary relation. Relations may depend on more variables. Constants and functions may also be used inside arguments of relations, and the functions are permitted to depend on an arbitrary number of variables. Thus one can make substitutions to create formulas like $\forall x \exists y\, G(\phi(x,y), \psi(\alpha(x), c))$, where $G$ is a binary relation, $\phi$, $\psi$ are functions of two variables, $x$, $y$ are variables, and $c$ is a constant symbol. Expressions like $\phi(x,y)$, $\psi(\alpha(x), c)$ are called *terms*, and in general terms may be constructed from constants and variables using any combination of functions. For arithmetic one permits additional symbols to represent numbers and basic arithmetic operations, and one adds axioms about the arithmetic objects. For example, in arithmetic $x + y$ is a term, $x + y = z$ is a formula, and there is an axiom to the effect that $x + 0 = x$.

The provable statements in arithmetic are called *theorems*. The same term is used in any context in which there are special axioms describing a mathematical structure, like arithmetic. Statements which are true independently of particular mathematical structure, as in ordinary propositional or predicate logic, are called *tautologies*.

Propositional and predicate logic are *sound* and *complete*. Roughly speaking this means that provability is equivalent to being true in all interpretations. In contexts with extra mathematical structure (like arithmetic) one should restrict oneself to interpretations which are compatible with the given structure.

These three systems are very different in terms of the levels of complexity that arise naturally in them. For this discussion we need to have some notions from complexity theory, but rather than descend into details let us just say a few words. The reader probably has already a reasonable idea of what is an algorithm. To make it more formal one has to be precise about what are the acceptable inputs and outputs. A typical output might simply be an answer of "YES" or "NO". The input should be encoded as a string of symbols, a "word" in the language generated by some alphabet. For instance, formulas in logic can be encoded in this way. A basic method to measure the complexity of an algorithm is to ask how long it takes the algorithm to give an answer when the input has length $n$ (= the number of symbols). Is the amount of time bounded by a polynomial in $n$, an exponential in $n$, some tower of exponentials, etc.

In propositional logic the problems that typically arise are resolvable by obvious algorithms in exponential time, and the difficult questions concern the existence of polynomial time procedures. For instance, given a propositional formula $A$, let us ask whether it is *not* a tautology. This question can be resolved in exponential time, by checking truth tables. This is actually an "NP" problem, which amounts to the fact that if $A$ is not a tautology, then there is a fast reason for it, namely a set of truth values for the propositional variables inside $A$ for which the result is "False". The difficulty is that this choice of truth values depends on $A$ and a priori one has to search through a tree to find it. This problem turns out to be NP-complete, which means that if one can find a polynomial-time algorithm for resolving it, then one can do the same for all other NP problems, such as the travelling salesman problem, or determining whether a graph has a Hamiltonian cycle.

By contrast in predicate logic one typically faces issues of algorithmic decidability or undecidability. That is, whether there is an algorithm at all that always gives the right answer, never mind how long it takes. The problem of determining whether a formula in predicate logic is a tautology is algorithmically undecidable. One can think of this as a matter of complexity, as follows. The tautologies in predicate logic of length at most $n$ is a finite set for which one can choose a finite set of proofs. Let $f(n)$ denote the maximum of the lengths of the shortest proofs of tautologies of length $\leq n$. The fact that there is no algorithm for determining whether or not a formula is a tautology means that $f(n)$ grows very fast, faster than any recursive function.

In order to find a proof of a tautology one has to permit a huge expansion. This is reminiscent of the word problem for finitely presented groups. Again there is no algorithm to determine whether a given word is trivial or not. One can think of this in terms of the huge (nonrecursive) expansion of a word which may be needed to establish its triviality.

This is a nice feature of proof theory and complexity. If one looks at the standard examples of complexity problems that deal with exponential versus polynomial time, they are usually very different from the standard examples of problems that are algorithmically decidable or undecidable. In proofs they can fit under the same roof in a nice way. "Herbrand's theorem" – discussed in Section 7 below – permits one to code predicate logic in terms of propositional logic, but with large expansion.

What about arithmetic? In arithmetic infinite processes can occur, coming from mathematical induction. One can still study the level of "infinite complexity" though. (This terminology may seem strange, but indeed the idea of measuring infinite complexity is present in many areas of mathematics, even if it is not always expressed as such. A lot of classical analysis can be seen in this way. One cannot describe a general function with a finite number of parameters, but one can do this approximately in the presence of some smoothness conditions, with the degree of approximation related to the degree of smoothness. With little or no smoothness even the approximate behavior must be counted in more infinite ways. The differentiability almost everywhere of Lipschitz functions provides a good example of this phenomenon. One encounters the necessity of infinite processes in topology as well.)

There are dynamics inside proofs, and these dynamics seem to be closely related to the matters of complexity just described, in all three cases, of propositional logic and predicate logic and arithmetic. Our next task is to describe a precise logical system in which we can work, and then explain the cut-elimination theorem, which provides an effective procedure for eliminating dynamics from proofs. We shall discuss some of the combinatorial aspects of cut-elimination, and some of its consequences. We shall give examples of proofs which use cuts in an interesting way, starting with the John-Nirenberg theorem in real analysis. The proof constructs an expanding tree of intervals in the real line by iterating a process that is coded in a single lemma. The proof has an interesting dynamical structure, in which an interval that is produced by the lemma is then fed back into the lemma to make more intervals, and so forth. This is an example of the kind of substitution that one expects in complicated proofs in predicate logic. We give a more elementary example with similar dynamical structure in Section 9. This example comes from [12] and is easier to formalize precisely. In Section 10 we explain the notion of the *logical flow graph* which provides a tool for seeing dynamical structure within proofs more clearly.

We would like to thank M. Baaz, M. Gromov, R. Kaye and C. Tomei for their comments and suggestions.

# 2 Sequent calculus

In our discussion of logic so far we have talked about the language – what is a formula? – but not about proofs. The choice of proof system is a nontrivial matter. We shall work with sequent calculus and the logical system $LK$. This turns out to have nice combinatorial properties, as we shall see.

The concept of a sequent can be a little confusing at first. A sequent is something of the form

$$A_1, A_2, \ldots, A_m \to B_1, B_2, \ldots, B_n$$

where the $A_i$'s and the $B_j$'s are formulas. The interpretation of this sequent is "from $A_1$ and $A_2$ and ... and $A_m$ follows $B_1$ or $B_2$ or ... or $B_n$". However, the arrow $\to$ is *not* a connective, nor is the sequent a formula.

Because we are using $\to$ for sequents in this manner, we use the symbol $\supset$ for the connective that represents implication. If $A$ and $B$ are formulas, then $A \supset B$ is also a formula, which is interpreted as "$A$ implies $B$".

So then what is the point of sequents? Why do we not simply use the formula

$$A_1 \wedge A_2 \wedge \ldots \wedge A_m \supset B_1 \vee B_2 \vee \ldots \vee B_n?$$

The two have the same interpretation, but they are different combinatorially. Even in ordinary reasoning we would not normally take all of our information $A_1$, ..., $A_m$ and package it into the single unit $A_1 \wedge \ldots \wedge A_m$. The commas permit us to treat the formulas $A_i$, $B_j$ as individuals which can each be used separately. We shall see that this flexibility is important in the system $LK$, for which we have certain *monotonicity* properties in the way that formulas are constructed.

We should say that in a sequent as above, the formulas $A_i$ and $B_j$ are permitted to have repetitions, and this turns out to be important. We do not care about the ordering, however. We might say "multisets" of formulas to make clear that we mean unordered collections in which repetitions are counted. We also permit empty sets of formulas, e.g.,

$$A_1, A_2, \ldots, A_m \to \qquad \text{and} \qquad \to B_1, B_2, \ldots, B_n$$

are permissible sequents. As a matter of notation we shall typically use upper-case roman letters $A, B, C...$ to denote formulas, and upper-case greek letters like $\Gamma, \Delta, \Lambda...$ to denote collections of formulas.

For a nice example of a sequent, consider $\Gamma \to \Delta$, where $\Gamma$ is the collection of formulas

$$\Gamma = \{\bigvee_{j=1}^{n} p_{i,j} : i = 1, \ldots, n+1\}$$

and $\Delta$ is the collection of formulas

$$\Delta = \{p_{i,j} \wedge p_{m,j} : 1 \leq i < m \leq n+1, 1 \leq j \leq n\}.$$

The $p_{i,j}$'s here are propositional variables. The sequent $\Gamma \to \Delta$ represents a finite version of the *pigeon-hole principle*: if you have $n+1$ balls, and each ball is placed in one of $n$ boxes, then there is at least one box which contains two balls. It is a valid sequent whose proof is somewhat tricky in propositional logic. Normally we think of it as a part of more powerful languages for which the proof is immediate, but in propositional logic it is more subtle, as we shall see.

Here now is the system $LK$. It should be interpreted as follows. Our basic objects are sequents, in the sense that we prove a sequent, we do not prove a formula. If we want to think of proving a formula $A$, then we should prove the sequent $\to A$. To prove a sequent we begin with axioms and derive new sequents from them using certain rules. For $LK$ the *axioms* are sequents of the form

$$A, \Gamma \to \Delta, A$$

where $A$ is any formula and $\Gamma, \Delta$ are any collections of formulas. The rules come in two groups, *logical rules* and *structural rules*. In these rules we write $\Gamma, \Gamma_1, \Gamma_2$, etc., for collections of formulas, and we write $\Gamma_{1,2}$ as a shorthand for the combination of $\Gamma_1$ and $\Gamma_2$. Remember always that we allow repetitions, so that one must count multiplicities when combining collections of formulas, and that we do not care about the ordering of the formulas on either side of the sequent.

The *logical rules* are used to introduce connectives, and they are given as follows:

$\neg$ : *left*  $\quad \dfrac{\Gamma \to \Delta, A}{\neg A, \Gamma \to \Delta}$  $\qquad\qquad$  $\neg$ : *right*  $\quad \dfrac{A, \Gamma \to \Delta}{\Gamma \to \Delta, \neg A}$

$\wedge$ : *right*  $\quad \dfrac{\Gamma_1 \to \Delta_1, A \quad \Gamma_2 \to \Delta_2, B}{\Gamma_{1,2} \to \Delta_{1,2}, A \wedge B}$

$\wedge$ : *left*  $\quad \dfrac{A, B, \Gamma \to \Delta}{A \wedge B, \Gamma \to \Delta}$

$\vee$ : *left*  $\quad \dfrac{A, \Gamma_1 \to \Delta_1 \quad B, \Gamma_2 \to \Delta_2}{A \vee B, \Gamma_{1,2} \to \Delta_{1,2}}$

$\vee$ : *right*  $\quad \dfrac{\Gamma \to \Delta, A, B}{\Gamma \to \Delta, A \vee B}$

$\supset: left$ $\quad \dfrac{\Gamma_1 \to \Delta_1, A \quad B, \Gamma_2 \to \Delta_2}{A \supset B, \Gamma_{1,2} \to \Delta_{1,2}}$

$\supset: right$ $\quad \dfrac{A, \Gamma \to \Delta, B}{\Gamma \to \Delta, A \supset B}$

$\exists: left$ $\quad \dfrac{A(b), \Gamma \to \Delta}{(\exists x)A(x), \Gamma \to \Delta}$ $\qquad \exists: right \quad \dfrac{\Gamma \to \Delta, A(t)}{\Gamma \to \Delta, (\exists x)A(x)}$

$\forall: left$ $\quad \dfrac{A(t), \Gamma \to \Delta}{(\forall x)A(x), \Gamma \to \Delta}$ $\qquad \forall: right \quad \dfrac{\Gamma \to \Delta, A(b)}{\Gamma \to \Delta, (\forall x)A(x)}$

The *structural rules* do not involve connectives and are the following:

$Cut$ $\quad \dfrac{\Gamma_1 \to \Delta_1, A \quad A, \Gamma_2 \to \Delta_2}{\Gamma_{1,2} \to \Delta_{1,2}}$

$Contraction$ $\quad \dfrac{\Gamma \to \Delta, A, A}{\Gamma \to \Delta, A}$ $\qquad \dfrac{A, A, \Gamma \to \Delta}{A, \Gamma \to \Delta}$

This is the system $LK$ for ordinary predicate logic. For propositional logic one drops the rules with quantifiers, and formulas are merely boolean combinations of propositional variables (with no functions or substitutions).

One should be a little careful about the rules for the quantifiers. In $\exists: right$ and $\forall: left$, any term $t$ is allowed which does not include a variable which lies already within the scope of a quantifier in the given formula $A$. In $\exists: left$ and $\forall: right$ one has the "eigenvariable" $b$ which should not occur free in $\Gamma, \Delta$.

This presentation of logic is quite different from what one ordinarily sees in an undergraduate course. It is easy to check however that the axioms and rules make sense in terms of usual reasoning.

The rule that probably looks the strangest is the cut rule. One can think of it as an elaboration on Modus Ponens, and again one can check that it is compatible with standard reasoning.

Although $LK$ does represent classical logic, it does so with some unusual subtleties, some properties of *monotonicity* and *conservation*. *Formulas never get simpler.* The logical rules permit us to make new formulas by combining old ones, and to introduce connectives, but they do not permit us to simplify formulas. *Formulas never disappear, except in cuts.* The only other simplification allowed is contraction, in which a repetition is reduced. *Formulas do not appear suddenly.* Everything has to be constructed from the formulas in the axioms.

An important consequence of this is that the size of a proof is controlled by the axioms, unless one stupidly applies the negation rules over and over again for no reason. In a large proof one should typically expect a large *number* of axioms. If the proof is large but the number of axioms is much smaller, and one does not repeat the negation rules in a foolish manner, then it means that there are a lot of "weak formulas", i.e., the formulas coming from the $\Gamma$'s and $\Delta$'s in the axioms $A, \Gamma \to \Delta, A$. If this happens then either the proof relied on foolish contortions involving weak formulas, so that it could be shortened, or the proof has a small number of steps and the endsequent is full of formulas coming from weak formulas, which is not too interesting.

Sometimes a crucial ingredient for the size of a proof looks dumb in terms of reasoning. The contraction rule, for instance, looks pretty silly, but we shall see in the next sections that it matters much for complexity. (In *Linear logic* [24], an extension of classical logic, contractions are controlled in an explicit way through *exponential* operators.) The problem of whether or not a sequent has a short proof can be seen as a purely combinatorial issue, apart from reasoning.

## 3 Cut elimination

**Theorem 1** (Gentzen [22, 25, 57]) *Any proof in $LK$ can be effectively transformed into a proof which never uses the cut rule. This works for both propositional and predicate logic.*

There is a version of this for arithmetic but one has to allow infinite proofs, because of induction. One can analyze the level of infiniteness used, measured by countable ordinals, as in [51].

This is a gorgeous theorem. It says first that we can make proofs without Modus Ponens, which is a bit striking. It also has nice consequences, including the *subformula property*: given a proof of a sequent which does not use the cut rule, then every formula which appears in the proof also appears as a subformula of a formula in the final sequent. This follows from the monotonicity properties of $LK$.

There is a nice image associated to proofs in $LK$ and the way that formulas are constructed. One can code the combinatorics of the construction of a formula with a tree. Combining two formulas $A$ and $B$ into $A \wedge B$, for instance, corresponds to combining the trees associated to $A$ and $B$. As formulas progress through $LK$ they are constructed in this manner. The trees never disappear when there are no cuts.

We shall discuss the proof of the cut-elimination theorem in a moment, but let us begin with the following question: what is the price for eliminating cuts? Just think about the difficulty in ordinary mathematical activity of making proofs which satisfy the subformula property.

The price comes in the expansion. The elimination of cuts leads to a simplification in the dynamical structure of the proof at the cost of large expansion in the size. There are propositional tautologies for which cut-free proofs must be exponentially larger than proofs with cuts, and in predicate logic the expansion can be non-elementary. See [44, 45, 53, 54, 55, 58].

A nicely concrete example of this expansion is provided by the *pigeon-hole principle*. We saw in the preceding section how the pigeon-hole principle can be formalized as a sequence of propositional sequents, one for each positive integer $n$. It turns out that the pigeon-hole principle can be proved in $LK$ with a proof of polynomial size in $n$ if one allows cuts [7], while exponential size is required for proofs without cuts [31]. See [1, 2, 49, 4, 3] for related work. Analogous results hold for the propositional version of the Ramsey theorem (see [47]).

The idea of lengths of proofs is quite amusing from the perspective of reasoning. It suggests that there are some statements that are true that we cannot understand in practice because it would take too long, and that there are statements which we can understand if we permit ourselves to use cuts and not otherwise. Induction plays a similar role in arithmetic, and indeed the strength of induction required in arithmetic for various purposes has been much studied. (See [30].)

What is it about the cut rule that permits this compression of proofs? One can think that the size of the proof is being compressed even if the dynamical content remains the same in essence. The use of lemmas permits one to make repeated substitutions with the same coding. We shall see this more concretely in Sections 8 and 9. There are graphs associated to proofs which trace the flow of logical occurrences, and these graphs are approximately trees for cut-free proofs, but proofs with cuts can have *cycles*. We shall discuss this further in Section 10.

As in the introduction, the existence of short proofs for all propositional tautologies is equivalent to $NP = co - NP$. It is not clear exactly what should create an impassable obstruction to compression. People have looked at combinatorial principles stated as tautologies in the hope of showing that the proofs had to be long even if cuts were allowed, but no one has succeeded in doing this so far. We already saw that the pigeon-hole principle does have short proofs (polynomial size). Similarly there are finite versions of Ramsey theorems which have been coded as statements in propositional logic and for which there are short proofs. Unlike the pigeon-hole principle these short proofs have not been given explicitly though. See [47].

In considering issues of complexity one should not be distracted too much by the interpretation in terms of reasoning. It seems more natural to view these as combinatorial problems which are merely expressed in the language of logic. Indeed they seem to have a natural geometry to them, discussed in Section 10.

How can one try to prove the cut-elimination theorem? The first point to understand is that it is not simply a matter of expressing the cut rule in terms of the others. When one studies logic as an undergraduate one is frequently told that various languages are equivalent, or various proof systems are equivalent, by dint of rules that permit translations from one to the other. That is not what is happening here. One has to operate on an actual proof, and the argument is more global in nature.

Roughly speaking the argument works by systematically going through a proof and simplifying the cuts. One starts at the bottom, near the final sequent, and works upwards, trying to push the cuts up until one encounters axioms or other simple situations. For instance, if we get to a cut of the form

$$\frac{A \to A \quad A \to B}{A \to B}$$

we can eliminate it stupidly by throwing away the axiom. In general though we have to take into account the specific structure of the situation. Consider the following situation.

$$\frac{\Gamma \to A \quad \dfrac{\dfrac{A \to B \quad A \to C}{A, A \to B \wedge C}}{A \to B \wedge C}}{\Gamma \to B \wedge C}$$

We used here first the rule for introducing a conjunction on the right, then a contraction, then a cut. Think of this as appearing in the middle of a proof, so that we have a proof $\Pi_0$ of $\Gamma \to A$, a proof $\Pi_1$ of $A \to B$, and a proof $\Pi_2$ of $A \to C$, and we are combining these three proofs to get a proof of $\Gamma \to B \wedge C$. We have used here a cut and we want to get rid of it. To do this we work as follows.

$$\frac{\dfrac{\Gamma \to A \quad A \to B}{\Gamma \to B} \quad \dfrac{\Gamma \to A \quad A \to C}{\Gamma \to C}}{\dfrac{\Gamma, \Gamma \to B \wedge C}{\Gamma \to B \wedge C}}$$

Here we started with a pair of cuts, followed by the conjunction rule and then contractions. This gives us another way to combine the proofs $\Pi_0$, $\Pi_1$, and $\Pi_2$ into a proof of $\Gamma \to B \wedge C$. From the point of view of cut elimination this arrangement is better than the previous one. That might seem strange, since we now have two cuts instead of one, but the two new cuts are at simpler stages in the proof than the original one. We have moved the cut up higher in the proof, across the contraction (from $A, A$ to $A$). In doing this we have introduced new contractions (to get from $\Gamma, \Gamma$ to $\Gamma$), but more importantly we have had to use the proof $\Pi_0$ of $\Gamma \to A$ twice. This is the principal source of *expansion* in the process of cut elimination, the duplication of subproofs that occurs when one pushes the cut up over a contraction.

This example corresponds to the idea of a *lemma* in ordinary mathematics. In the original piece of proof we knew that $A$ would imply each of $B$ and $C$, we had a derivation of $A$ from $\Gamma$, and we wanted to obtain $B \wedge C$. We did not want to have to derive $A$ from $\Gamma$ twice, we wanted one "lemma" to say that it was true. By using the contraction and cut we were able to do this, but by eliminating the cut we had to repeat the proof for each application.

The use of the cut enabled us to make a shorter proof, but to gain this efficiency we had to merge two occurrences of $A$ even though they could be used in completely different ways. To understand this point it is helpful to think of arithmetic. Think of $A$ as saying that some property of numbers is preserved when one multiplies two of them together. This fact might be employed several times in the proof, applied to different numbers or terms, even though proved only once. It may be that one lemma is applied to terms obtained from previous applications of itself. This type of phenomenon occurs in the examples discussed in Sections 8 and 9.

Here is another example, but in predicate logic this time.

$$\cfrac{\cfrac{\cfrac{\cfrac{A \to F(t)}{A \to \exists x F(x)} \quad \cfrac{B \to F(s)}{B \to \exists x F(x)}}{A \vee B \to \exists x F(x), \exists x F(x)}}{A \vee B \to \exists x F(x)} \quad \cfrac{F(a) \to \Delta}{\exists x F(x) \to \Delta}}{A \vee B \to \Delta}$$

Again think of this as being part of a larger proof, in which we have already proofs of $A \to F(t)$, $B \to F(s)$, and $F(a) \to \Delta$, where the *eigenvariable* $a$ does not occur free in $\Delta$. Here $s$ and $t$ are *terms*, which means that they can be constructed from constants and variables using function symbols. This provides a nice example of the way that the cut rule can be used; it permits us to make the lemma that $\exists x F(x) \to \Delta$ rather than deriving $\Delta$ from $F(x)$ for each possible choice of $x$. To push the cut upward we need to make separate proofs for $s$ and $t$. We do this as follows.

$$\cfrac{\cfrac{\cfrac{A \to F(s) \quad F(s) \to \Delta}{A \to \Delta} \quad \cfrac{B \to F(t) \quad F(t) \to \Delta}{B \to \Delta}}{A \vee B \to \Delta, \Delta}}{A \vee B \to \Delta}$$

We get the proofs of $F(s) \to \Delta$ and $F(t) \to \Delta$ by taking the proof of $F(a) \to \Delta$ in the original and substituting the terms $s$ and $t$ for the eigenvariable $a$.

This ability to make substitutions is a fundamental difference between predicate and propositional logic. Substitutions interact with the cut rule in a very substantial way. In passing from $F(s)$ and $F(t)$ to $\exists x F(x)$ we are doing something very strong, since the terms $s$ and $t$ need not have anything to do with each other.

These two examples illustrate the ideas and effects of pushing a cut up across a contraction. In fact there is a general procedure which works in all cases. If we start with

$$\cfrac{\cfrac{\overset{\Pi_1}{\Gamma_1 \to \Delta_1, A, A}}{\Gamma_1 \to \Delta_1, A} \quad \overset{\Pi_2}{A, \Gamma_2 \to \Delta_2}}{\Gamma_1, \Gamma_2 \to \Delta_1, \Delta_2}$$

where $A$ is an arbitrary formula, then we can replace it with

$$\cfrac{\cfrac{\overset{\Pi_1}{\Gamma_1 \to \Delta_1, A, A} \quad \overset{\Pi_2}{A, \Gamma_2 \to \Delta_2}}{\Gamma_1, \Gamma_2 \to \Delta_1, \Delta_2, A} \quad \overset{\Pi_2}{A, \Gamma_2 \to \Delta_2}}{\cfrac{\Gamma_1, \Gamma_2, \Gamma_2 \to \Delta_1, \Delta_2, \Delta_2}{\vdots \; contractions \atop \Gamma_1, \Gamma_2 \to \Delta_1, \Delta_2}}$$

In doing this we have to duplicate the proof $\Pi_2$. This can lead to large expansion in the size of the proof if we have to do it many times.

This explains how one can push a cut up across any contraction. Consider the problem of pushing a cut up over a conjunction, as in

$$\cfrac{\cfrac{\Gamma_1 \to \Delta_1 A \quad \Gamma_2 \to \Delta_2 B}{\Gamma_1, \Gamma_2 \to \Delta_1, \Delta_2, A \wedge B} \quad \cfrac{A, B, \Gamma \to \Delta}{A \wedge B, \Gamma \to \Delta}}{\Gamma, \Gamma_1, \Gamma_2 \to \Delta, \Delta_1, \Delta_2}$$

We can replace this with

$$\cfrac{\Gamma_2 \to \Delta_2, B \quad \cfrac{\Gamma_1 \to \Delta_1, A \quad A, B, \Gamma \to \Delta}{B, \Gamma, \Gamma_1 \to \Delta, \Delta_1}}{\Gamma, \Gamma_1, \Gamma_2 \to \Delta, \Delta_1, \Delta_2}$$

Now we are using two cuts, but they have lower complexity.

For this replacement we could just as well have used

$$\cfrac{\Gamma_1 \to \Delta_1, A \quad \cfrac{\Gamma_2 \to \Delta_2, B \quad A, B, \Gamma \to \Delta}{A, \Gamma, \Gamma_2 \to \Delta, \Delta_2}}{\Gamma, \Gamma_1, \Gamma_2 \to \Delta, \Delta_1, \Delta_2}$$

This reflects an important feature of cut-elimination: there is no canonical way to do it. We had a choice here about how to arrange the cuts. Similarly for contractions, if both appearances of the cut formula $A$ were obtained from contractions, then we would have a choice as to which subproof to duplicate first. In principle we can have procedures of cut-elimination which go on forever. Of course the point of the theorem is that one can always find a way to eliminate cuts in a finite number of steps. One can even make deterministic procedures by imposing conditions on the manner in which the transformations are carried out. See [24].

The cases that we have considered illustrate well the general scheme of the proof of the cut-elimination theorem. The principle is to push the cuts up higher in the proof, but we have to be careful about the notion of "progress", because we typically increase the number of cuts at each stage of the process. In the examples about contractions we made progress in the sense that we reduced the number of contractions above the cut formula, even though we may increase the total number of contractions by adding them below the cut. In the example with conjunctions we reduced the complexity of the cut formula. It is not hard to make examples to exhaust the possibilities, but a complete proof requires a tedious verification of cases that we shall not provide. (See [25, 57].)

# 4 Mathematics and formal proofs

Mathematical logic provides a way to formalize ordinary mathematical activity, and cut-elimination has a very interesting role in this.

For this story we should back up to the *Hilbert program*, which sought to show that mathematics could be formalized in a pure way and that in principle abstraction could be avoided. (See [6].) One can try to treat mathematics in a completely symbolic manner, where mathematical formulas are strings of symbols constructed through fixed formal rules. The rules of logical inference are formal as well and permit the passage from one string of symbols to another. In this way proofs can be treated as mathematical objects in their own right and studied mathematically.

In mathematics we have a natural informal hierarchy of abstraction. We can start with the very concrete, like natural numbers and finite graphs, and proceed to the more transcendental, real and complex numbers, real and complex analysis. We can proceed further from "simple" infinite processes to infinite-dimensional ones, like Hilbert spaces and operators on them, then further still to highly abstract and nonconstructive concepts, often associated to the axiom of choice, such as ultrafilters. Many aspects of analysis, for instance, are touched by the phenomenon of nonconstructive existence through compactness or the Hahn-Banach theorem.

In some cases these abstractions are "approximately concrete", as with infinite processes which admit well-behaved finite approximations. Sometimes remote abstractions wander into more concrete worlds. Compactness or the Hahn-Banach theorem may be used to "provide" theoretical solutions to explicit differential equations, or compactness can lead to uniform bounds without providing a clue as to how to generate an actual number. Even more elementary questions about integers are sometimes treated through transcendental methods or nonconstructive abstractions of compactness.

These phenomena are naturally troubling. The literature is full of debates and attempts to address the problem. Hilbert's program, in its strongest form, would have provided a very attractive resolution of these difficulties. It sought to show, roughly speaking, that concrete statements that are true should always have finite proofs. That the infinite methods of general mathematical activity would not lead us to trouble, that the infinite abstractions were a convenience that could, in principle, be replaced with more direct elementary arguments for elementary assertions.

In its strongest form Hilbert's program failed and failed utterly, because of the celebrated work of Gödel.

In some ways Gödel's work has had unfortunate negative side effects. Gödel's results have enormous conceptual consequences, and they force one to confront and accept some troubling phenomena, but some of the ideas behind Hilbert's program retain their strength despite being overshadowed by the failure of other aspects.

Gentzen's theorem arose in this context. One can view it as providing a positive result in contrast with Gödel's work. Cut elimination gives an approach to converting indirect proofs into direct ones. Gentzen's work also helps to clarify the precise meaning of an "elementary" proof, which Hilbert had left vague and intuitive.

Cut elimination need not work in an arbitrary logical system, or it may work with qualifications. In arithmetic it may convert a finite proof into an infinite one. Still, one can often control the level of infinite processes used (transfinite induction up to a certain ordinal) and thereby obtain "quasi-elementary" proofs of elementary statements. One of Gentzen's purposes was to bring the consistency of arithmetic closer to elementary arguments.

Keep in mind that cut elimination does transform finite proofs into finite ones in ordinary predicate logic. Arithmetic is simply much more complicated, as we see in Gödel's results and the necessity of passing to infinite proofs when eliminating cuts. Remember though that the set of tautologies in ordinary predicate logic is algorithmically undecidable, and that cut elimination leads to *large* expansion of proofs. Thus ordinary predicate logic is finite in ways that arithmetic is not, but the expansion is still there, in a weaker form.

This view of cut elimination and Hilbert's program is illustrated well by some work of Girard [25]. The story begins with Hilbert's program in reverse: Furstenberg and Weiss [20] found a very elegant proof through transcendental methods of dynamical systems of the van der Waerden theorem [59] on arithmetic progressions. In this case the elementary proof came first, the short nonelementary proof arrived much later. Girard showed however that one could recover the elementary combinatorial arguments by applying the procedure of cut elimination to the methods of Furstenberg and Weiss.

Other examples of analysis of mathematical proofs through cut elimination can be found in [5, 32, 33, 34] and in unpublished work of Kreisel concerning a theorem of Littlewood in number theory. For this type of analysis a basic tool is the *no counterexample interpretation* of Kreisel for turning infinite arguments into finite ones through the introduction of functionals.

This basic idea, of looking for consequences of Gentzen's work for ordinary mathematics, occurs repeatedly in the writings of Kreisel [37, 39, 38], and he raises many concrete questions. Progress has been made, but the issue remains to be addressed in a strong way.

# 5    Remarks about the subformula property

A proof satisfies the *subformula property* if every formula that appears in the proof also arises as a subformula of a formula in the endsequent. We saw in Section 3 that any provable sequent has a proof which enjoys the subformula property, because any cut-free proof has this property.

We should be a little bit careful here. What exactly is a subformula? Basically $A$ is a subformula of $B$ if you can get $B$ from $A$ by adding things to it. For propositional logic this is straightforward, but quantifiers in predicate logic bring a subtlety with them. We can start with a formula $R(t)$, where $F$ is a unary relation and $t$ is a term, and we can build from it the formula $\exists x R(x)$. We consider $R(t)$ to be a subformula of $\exists x R(x)$. If $s$ is another term then $R(s)$ is also considered to be a subformula of $\exists x R(x)$.

This is an enormous difference between propositional and predicate logic. In propositional logic if we know a formula then we know its ancestors exactly. If it were true in predicate logic that the ancestors of a formula are determined by the formula itself, then there would be an algorithm for deciding which sentences are tautologies, which is not true. This would amount to saying that any tautology would have a cut-free proof which could then be controlled.

The size of the formulas which appear in a cut-free proof in propositional logic is bounded by the size of formulas in the endsequent. In predicate logic this does not work because the terms that appear inside the proof can be larger than the ones in the endsequent. In both propositional and predicate logic the *size* of a cut-free proof – the number of symbols in the proof – can be controlled in terms of the size of the endsequent together with the number of lines in the proof. (See [35].) In both cases the number of lines in a cut-free proof can be very large compared to the size of the endsequent because of the presence of contractions.

We should emphasize the relation between the subformula property and the idea of "lemmas" in a proof. Gentzen's cut-elimination theorem is often described as a procedure for making *direct* proofs, for avoiding intermediate results which are more general than the final theorem. The subformula property is a manifestation of this idea, that nothing occurs in the proof which is more general than the final result. To understand this in a more concrete way think about proving a formula $A(t)$ that has no quantifiers but does have a term $t$ in it. In establishing this tautology we might use formulas of the form $R(s_i)$ many times, for various terms $s_i$. One can imagine that a shorter proof might be possible using $\forall x R(x)$. This would not be allowed in a cut-free proof, because of the subformula property, and because $A(t)$ has no quantifiers. In a cut-free proof all occurrences of $R(s_i)$ for the relevant terms $s_i$ would have to be listed separately.

Some of these phenomena occur in the examples discussed in Sections 8 and 9 below.

# 6 The Craig interpolation theorem

Roughly speaking the Craig interpolation theorem [17] states that if one can prove the sequent
$$A \to B$$

in $LK$ then one can also prove the sequents

$$A \to C \quad \text{and} \quad C \to B$$

where the formula $C$ – called the *interpolant* of $A$ and $B$ – involves only the language common to $A$ and $B$. *Common language* means the common propositional variables in the case of propositional logic, and the relations, functions, and constants for the predicate case.

In terms of reasoning this is not at all surprising. If $A$ involves apples and oranges, $B$ involves apples and bananas, and $A$ implies $B$, then $A$ ought to imply a statement that involves only apples, and $B$ ought to follow from a statement that involves only apples. The oranges should not help, and the bananas should not hurt.

So what is the mystery then? The Craig theorem is trickier to prove than one might think. One has to have the same statement about apples for both $A$ and $B$! To construct $C$ and proofs of $A \to C$ and $C \to B$ the cut-elimination theorem is extremely useful. Once one has a cut-free proof, the construction of $C$ is a fairly simple combinatorial problem (see [40] and also [25, 57]).

In fact one can formulate Craig interpolation in purely combinatorial terms [10], in such a way that the special nature of formulas does not really matter. The interpolation theorem can be seen as a general statement with few structural requirements that would apply as well to graphs or polygons or tessellated surfaces as to formulas in logic.

It is important here that we have cut elimination to start with. It would be much more difficult to give a general combinatorial formulation of Craig interpolation for proofs with cuts. With cut-elimination as a starting point one need only understand how to combine basic objects in a natural way. Formulas are the basic objects in logic, but to have a more combinatorial image it is nicer to think of them simply as trees. Think of coding the construction of a formula by a tree, and then forgetting the logic and just keeping the tree. It is easy to combine two trees, by connecting them at the root, and after cut elimination this is essentially the kind of operation that Craig interpolation requires. If there are cuts in the proof then the matter is entirely different. One would then have to go inside the structure of the objects, rather than simply combining them.

This brings us to an interesting point about cut-elimination, as an indication of simplicity or "finiteness" of the combinatorics of proofs. It is better to have cut elimination with large expansion than to not have it at all. In what other contexts in mathematics is there something analogous? Some "normalization" which always exists (in a nontrivial way)? One is tempted to draw an analogy with singularities of algebraic varieties. Resolution of singularities is a difficult matter, but it reflects an underlying finiteness that is not typically present in analysis or topology, for instance.

Cut elimination permits us to derive some general results through effective constructions, like Craig Interpolation and Herbrand's theorem (discussed below). However, the use of cut elimination does not provide interesting bounds. One constructs the interpolant by induction, at each step combining old interpolating formulas to get a new one. Contractions in the proof do not lead to contractions of interpolating formulas. In fact the interpolant reflects the entire construction in the proof. For this reason its size may be very large compared to the size of the endsequent, but linear in the size of the cut-free proof.

It is not apparent that interpolation should necessarily be as complicated as cut elimination.

What kind of bounds on the size of the interpolant $C$ can one obtain in terms of the sizes of $A$ and $B$? If the interpolant $C$ must be large compared to $A$ and $B$, what would that imply about the structure of $A \to B$?

In first-order predicate logic (with equality) one can have nonrecursive expansion for the size of the interpolant over the size of the endsequent, as in [19, 41].

What about propositional logic? Is it true that the size of $C$ is always bounded by a polynomial in the sizes of $A$ and $B$? If so, then one would have a general result in complexity theory, namely $NP \cap co-NP \subseteq P/poly$. This problem is a weaker cousin of $P = NP$, but it is equally unknown. The main point is that a language which is in $NP \cap co - NP$ can be coded in terms of a family of sequents $A_n \to B_n$ of controlled size, $n = 1, 2, 3, \ldots$, and a family of interpolants $C_n$ involving only the common variables of $A_n$ and $B_n$ would lead to a family of circuits which characterize the original language. See [42, 43, 11] for more details.

If one does not believe in polynomial bounds for general classes of problems which appear to require exponential time, then one would expect the failure of polynomial bounds for Craig interpolation.

This brings us back to our earlier questions about what is a difficult proof, what are the mechanisms by which proofs with cuts can be much shorter than cut-free proofs, etc. These questions fit well with the problem of knowing what structure in a sequent $A \to B$ is needed if the interpolant $C$ is always large. See [11] for some results in this direction.

What kinds of bounds can one get for the sizes of the proofs of $A \to C$ and $C \to B$ in terms of the proof of $\Pi : A \to B$? In the propositional case, for instance, can one construct $C$ and proofs of $A \to C$ and $C \to B$ whose sizes are bounded by a polynomial in the size of $\Pi$? Is there a polynomial-time algorithm for finding these proofs? These are all questions that are open. See [11] for further discussion of these matters.

There is an interesting variant of these questions, in which one starts with a truth assignment $\sigma$ for the variables common to $A$ and $B$ and looks for either a proof of $A^\sigma \rightarrow$ or $\rightarrow B^\sigma$ of controlled size. Roughly speaking $A^\sigma$ and $B^\sigma$ are obtained from $A$ and $B$ using the truth assignment $\sigma$. See [11] for details. The idea is that if we have an interpolating formula $C$, then $\sigma$ converts $C$ simply into "true" or "false", and proofs of $A \rightarrow C$ and $C \rightarrow B$ would then give rise to either a proof of $A^\sigma \rightarrow$ or of $\rightarrow B^\sigma$ according to whether we obtained "false" or "true" from $C$, respectively. If for each truth assignment $\sigma$ we can decide which of $A^\sigma \rightarrow$ or $\rightarrow B^\sigma$ is provable in polynomial time, then that is essentially the same as finding an interpolant $C$ of polynomial size. In both cases we get a function which can be computed in polynomial time even if the descriptions are not the same.

# 7 Herbrand's theorem

Let us begin with a simple version of Herbrand's theorem.

**Theorem 2** *If $A(x)$ is a formula without quantifiers in which the variable $x$ appears free, then $\rightarrow \exists x A(x)$ is provable in LK if and only if there is a finite collection of terms $t_1, \ldots, t_n$ such that $\rightarrow A(t_1), \ldots, A(t_n)$ is provable using only propositional rules (i.e., without quantifier rules).*

We should emphasize that $A(x)$ is allowed to be a *formula*, and not simply a relation, so that $A(x)$ might have the form $F(x) \wedge \neg G(x, \psi(x))$, for instance.

Roughly speaking, the theorem says that if you can prove $\exists x A(x)$, then you can make the existence explicit by producing a *finite* collection of terms, at least one of which satisfies the property $A(\cdot)$. Of course the converse is true because of the $\exists : right$ and contraction rules.

The theorem is very easy to obtain once we have cut elimination. If $\rightarrow \exists x A(x)$ is provable, then it is provable without cuts. Thus we may assume that we have a cut-free proof $\Pi$ of $\rightarrow \exists x A(x)$ to begin with. In particular $\Pi$ enjoys the subformula property. This means that any formula that occurs in $\Pi$ and which has a quantifier in it must be an occurrence of $\exists x A(x)$.

Let us start from the bottom of $\Pi$, at the endsequent $\rightarrow \exists x A(x)$, and think about what happens as we go up in the proof. We may pass through contraction rules, which would cause $\exists x A(x)$ to be duplicated as we go up. The first moment at which $\exists x A(x)$ is derived from something other than a contraction rule it must be obtained using $\exists : right$, starting from a formula of the form $A(t)$ for some term $t$.

It is easy to use these observations to reorganize $\Pi$ to get a proof of $\rightarrow A(t_1), \ldots, A(t_n)$ as desired, where $t_1, \ldots, t_n$ are the terms that arise when we peel off the existential quantifiers. More formally, one can first remove the contractions "at the bottom" and make some minor reorganizations as needed to get a proof of $\rightarrow \exists x A(x), \ldots, \exists x A(x)$, where none of the occurrences of $\exists x A(x)$ was obtained using a contraction, and then one peels off the quantifiers to get $\rightarrow A(t_1), \ldots, A(t_n)$. Notice that only propositional rules could have been used to obtain the $A(t_i)$'s, because of the subformula property (since $A(x)$ is quantifier-free).

The basic principles of this argument are very general. Let us work with general formulas in "prenex" form, which means that all the quantifiers are on the outside. (Thus $\exists x \forall y (F(x) \wedge G(y))$ is prenex, $(\exists x F(x)) \wedge (\forall y G(y))$ is not. Formulas can always be put into prenex form.) Suppose that a sequent $\Gamma \rightarrow \Delta$ consists only of prenex formulas and has a cut-free proof $\Pi$ in $LK$. The *midsequent theorem* [22, 23, 25] states that we can modify $\Pi$ to get a new proof $\Pi'$ which has the property that as soon as one of the quantifier rules is used in the proof, only quantifier rules and contraction rules can be used afterwards. One should think of a proof as being coded by a tree here, with different branches being independent of each other until the moment that they meet on the way to the endsequent. It is easy to see how the midsequent theorem is established. Because a cut-free proof enjoys the subformula property, all formulas that appear in it are in prenex form. Once a formula has a quantifier attached to it, one cannot use it actively in a logical rule that is not a quantifier rule. This makes it easy to rearrange the proof in order to use all the quantifier rules after the nonquantifier rules.

The midsequent theorem shows that in principle the quantifier rules are not essential for predicate logic. In fact there is a general statement to the effect that tautologies in predicate logic can be converted into propositional tautologies without losing information. We gave a version of this in Theorem 2, but it is not as simple to accommodate alternation of universal and existential quantifiers. To handle the general case one uses function symbols which were not in the original language, called *Skolem functions*. These special functions code relationships between terms in the Herbrand disjunction which allow one to recover the statement in predicate language. See [25] for more information.

Let us illustrate the meaning of Herbrand's theorem in ordinary mathematics with the following scenario. Suppose that we are interested in something like the word problem in finitely presented groups (or semigroups, which are technically simpler for this discussion). We might then be interested in

a sequent of the form
$$\forall x_1 R_1(x_1), \ldots, \forall x_k R_k(x_k) \to B(t)$$
where the formulas $\forall x_i R_i(x_i)$ contain information from the relations of the group (in the form that an element of the group multiplied by a relation gives back the same element of the group), and where $B(t)$ is a formula without quantifiers which contains the information that a given word, represented here as a term $t$, is trivial in the group.

This sequent is in almost the same form as for the above formulation of Herbrand's theorem. That is, we can use the negation rules to convert the universal quantifiers on the left side into existential quantifiers on the right side. We have several formulas instead of just one, but they can be combined on the right with disjunctions. We have several quantifiers, but that does not really matter.

If one can prove such a sequent, then it means that one can show that the word is trivial in the group using the given relations. Of course if one can prove this one should not use the relations more than a finite number of times, even though formulas of the form $\forall x_i R_i(x_i)$ allow for infinitely many possible choices of $x_i$ a priori. Herbrand's theorem is a general statement to this effect. Cut-elimination provides a procedure to make explicit the relations that are needed. Proofs with cuts can be shorter by not making explicit the way that the relations are used, or how often. This is the power of quantifiers; they are "infinite" objects which can enable one obtain "finite" conclusions much more quickly.

We shall look more closely at an example of this type in Section 9, in connection with the problem of defining large numbers in arithmetic. Before we do that we discuss another example from analysis.

# 8 The John-Nirenberg theorem

The John-Nirenberg theorem is a result in real analysis whose proof provides an interesting example for cut elimination and the combinatorics of proofs.

Let $f(x)$ be a locally integrable function on the real line. Given an interval $I$ in $\mathbf{R}$, write $A(f, I)$ for the average of $f$ over $I$,
$$A(f, I) = \frac{1}{|I|} \int_I f(y) \, dy,$$
where $|I|$ denotes the length of $I$. Define the "mean oscillation of $f$ over $I$" by
$$\mathrm{mo}(f, I) = \frac{1}{|I|} \int_I |f(x) - A(f, I)| \, dx.$$

We say that $f$ has *bounded mean oscillation* ($BMO$) if

$$\|f\|_* := \sup_I \mathrm{mo}(f, I) < \infty.$$

Bounded functions satisfy this property, but some unbounded functions do too, such as $f(x) = \log|x|$. Note however that $|x|^\alpha$ does not lie in $BMO$ for any $\alpha \neq 0$.

Think of the restriction of $f$ to an interval $I$ as being a localized snapshot of $f$. The $BMO$ condition provides a way to say that these snapshots all have bounded averages, except that we permit ourselves to subtract off the mean value as a kind of normalization.

How big can a $BMO$ function really be? A simple fact is that

$$|\{x \in I : |f(x) - A(f, I)| > \lambda\}| \leq \frac{\|f\|_*}{\lambda} |I|$$

for all $\lambda > 0$, where we write $|E|$ for the Lebesgue measure of $E$. This is just Tchebychev's inequality,

$$\lambda |\{x \in I : |f(x) - A(f, I)| > \lambda\}| \leq \int_I |f(x) - A(f, I)| \, dx$$

simply by definition. However, it turns out that for a $BMO$ function we actually have *exponential* decay of the measure as a function of $\lambda$, that is

$$|\{x \in I : |f(x) - A(f, I)| > \lambda\}| \leq 2^{-(\lambda-1)/4} \|f\|_* |I|$$

for all $\lambda > 0$. This means that $BMO$ functions are closer to being bounded than it might appear at first. This estimate is roughly the best possible, because the logarithm has exactly exponential decay on the interval $[0, 1]$, for instance.

This exponential decay is a famous result of John and Nirenberg. See [21, 56, 52]. It is very important that in the definition of $BMO$ we take the supremum over *all* intervals. We certainly do not get exponential decay on an interval $I$ simply from the knowledge that $\mathrm{mo}(f, I)$ is finite for that particular interval.

This is the general setting of the John-Nirenberg theorem. Let us now sketch the proof. The main lemma comes from a construction of Calderón and Zygmund, and it says the following. Suppose that we have an interval $J$ and a function $h$ on $J$ such that

$$\frac{1}{|J|} \int_J |h(y)| \, dy \leq 1.$$

Then we can find a subset $E$ of $J$ with the following three properties. First, $E$ has at least half the points of $J$,

$$|E| \geq \frac{1}{2}|J|.$$

Second, $h$ is not too big on $E$,

$$|h(x)| \leq 2 \qquad \text{for almost all } x \in E.$$

The third property is more geometric. It says that $J \backslash E$ can be realized as the disjoint union of a collection of subintervals $\{J_i\}$ of $J$, in such a way that

$$\frac{1}{|J_i|} \int_{J_i} |h(y)| \, dy \leq 4$$

for each $i$.

To understand what all of this means think about the first two properties and forget about the third for a moment. To get these properties we do not need any kind of special geometry, it is simply a question of measure theory. If $|h| \leq 1$ on average, then we cannot have $|h| > \frac{1}{2}$ on more than half the points.

This measure-theoretic argument says nothing about what happens on the set where $|h| > \frac{1}{2}$. The point of the Calderón-Zygmund construction is to break up this "bad" set into pieces where the average is bounded again, as stated in the third property. We cannot say exactly what happens on the bad set, but we can organize it in a good way.

Let us describe the main point of the Calderón-Zygmund construction. Start with $J$ and break it into its two halves. On each of these ask the question "Is the average of $|h|$ larger than 2?" When the answer is yes put that interval aside. It will become one of the $J_i$'s. When the answer is no take the given interval, cut it into halves, and ask the same question about the two new intervals.

Repeating this process indefinitely we get a bunch of intervals on which the average was larger than 2. Outside these intervals $|h| \leq 2$ a.e. The intervals are disjoint, and one can show that the sum of their measures is $\leq \frac{1}{2}|J|$, using the assumption that the average of $|h|$ on $J$ is $\leq 1$. The averages of $|h|$ on these "bad" intervals is $\leq 4$. This comes because we called an interval bad at the first moment that it was bad, it came from splitting an interval on which the average of $|h|$ was $\leq 2$.

This is roughly how the Calderón-Zygmund construction works. Let us explain how one proves the John-Nirenberg theorem. Fix an interval $I$, and assume for simplicity that $\|f\|_* \leq 1$. We apply the Calderón-Zygmund construction to $h = f - A(f, I)$ on $I$. We get a set on which $|f - A(f, I)| \leq \frac{1}{2}$ and a bunch of intervals $\{I_i\}$ on which the average of $|f - A(f, I)|$ is bounded by 4. The main point is to repeat the construction on each $I_i$, applied to the function $f - A(f, I_i)$. The assumption that $\|f\|_* \leq 1$ implies that the average of $|f - A(f, I_i)|$ over $I_i$ is bounded by 1. To relate what we get back to $f$ we use the fact that the average of $|f - A(f, I)|$ over each $I_i$ is bounded by 4, whence
$$|A(f, I_i) - A(f, I)| \leq 4.$$

We repeat the Calderón-Zygmund construction on each $I_i$, for each one we get a new family of intervals, we apply the construction to each of those, and repeat the process indefinitely. One can make computations to derive the exponential decay mentioned earlier.

For us the Calderón-Zygmund construction is the "lemma" that we are applying repeatedly. If we do not make cuts then we have to repeat the construction on each interval. Instead we can think of proving it once but using it many times, through contractions and cuts. In each application it is applied to different functions and intervals, and indeed it is applied to intervals that were produced in an earlier application of the lemma.

We shall not try to make a precise formalization of the proof in logic, but it should be clear in principle how the elimination of cuts corresponds to performing the Calderón-Zygmund construction explicitly each time it is needed, and that with contractions and a cut one can make a smaller proof with a lot of "coiling" or "cycling" in the proof, coming from the output of the lemma going back into it as input. This idea can be made more precise through the logical flow graph, which we discuss in Section 10.

If one wants to treat the proof of the John-Nirenberg theorem more precisely in terms of logic, then it might be better to work with finite versions of it. Instead of working with functions on the whole real line let us think of functions on the set of integers $1, 2, 3, \ldots, 2^N$, where $N$ is a large positive integer. We still have a natural notion of intervals in this case, although for technical reasons it is better to restrict ourselves to intervals for which the number of elements inside is an integer power of 2. This makes it easier to cut them in half. Instead of Lebesgue measure we use counting measure, so that the integrals become sums.

If we do this then we can get the same kind of inequalities as before. For the purposes of analysis the precise constants like 2 and 4 are not so important, what matters are uniform estimates, with constants which do not depend on $N$ in this discrete model.

In this discrete model one can formalize the John-Nirenberg inequality and its proof in fairly simple logical terms. One does not need anything like the full structure of the real line or measure theory, one could also restrict oneself to functions which take values in the rationals to avoid any use of the real numbers.

This discrete version of the John-Nirenberg theorem also contains all the information of the original. Notice that there is a small additional subtlety for the proof. In the real line all intervals are practically the same, one can change from one to another using an affine mapping. In the discrete model that is not quite the case, one cannot use dilations to identify intervals of different length. If one were to formalize a proof with cuts this might require different proofs of the Calderón-Zygmund construction for intervals of different lengths, rather than having one proof for all intervals. This could lead to a requirement for $N$ different lemmas in the proof instead of just one, one lemma for each dyadic length $2^j$.

This expansion in the number of lemmas reflects an interesting difference between "continuous" and "discrete" mathematics. There can be more symmetries in continuous mathematics which enable one to make shorter proofs.

Notice the strong role of quantifiers in the proof of the John-Nirenberg theorem. The argument works because we have information about *all* intervals, as in the condition $\|f\|_* \leq 1$. This is a common phenomenon in harmonic analysis, that properties of functions or sets are naturally expressed in terms of universal quantifiers which reflect crucial scale-invariance. In the proof we used universal quantifiers to avoid dealing with specific intervals. This lead to a kind of cycling in the argument. In the next section we describe a different example which is easier to analyze and which exhibits similar features of dynamics in short proofs with cuts.

## 9 Large numbers

We want to provide now a more precise example where quantifiers and the cut rule can be used to make short proofs, and where explicit proofs would need to be much longer. We do this in the context of arithmetic, through the concept of "feasible numbers" [46]. The logical flow graph (Section 10) associated to the proof described below has a rich structure of cycles.

To try to avoid overloading the reader with syntax let us not review the manner in which arithmetic is formalized in predicate logic. In short one specifies symbols to reflect the basic objects and operations in arithmetic $(0, =, <, +,$ and $*$ for multiplication), and one adds axioms to encode their basic properties. One of the main objects used in this formalization is the function $s(x)$, the *successor* function, which is interpreted in normal mathematics as $x + 1$. This is important for the formalization of mathematical induction, which we shall *not* use here.

We shall also "cheat" a little and add to our version of arithmetic the operation of exponentiation $z^x$ with the usual property that $(z^x)^y = z^{x*y}$. For the present purpose dealing with the exponential function at the level of first principles would be a distraction, but it is certainly possible to define it in the framework of Peano arithmetic and to establish its properties.

Let us add a new object to arithmetic, namely the symbol $F$ which represents a unary relation of *feasibility*. The idea is that if $x$ represents a natural number, then $F(x)$ represents the statement that $x$ can be constructed in some feasible manner. The properties of $F$ can be described as follows:

$$
\begin{array}{ll}
 & F(0) \\
F : equality & x = y \supset (F(x) \supset F(y)) \\
F : inequality & F(x) \wedge (y < x) \supset F(y) \\
F : successor & F(x) \supset F(s(x)) \\
F : plus & F(x) \wedge F(y) \supset F(x + y) \\
F : times & F(x) \wedge F(y) \supset F(x * y)
\end{array}
$$

To show that a number $n$ is feasible one ought to be able to construct a proof of $F(n)$. We do not allow induction to be used over formulas containing $F$. Otherwise we could prove $\forall x F(x)$ in a few steps.

Note that we are *not* allowing an exponential rule for $F$.

In mathematical logic one sometimes looks at these axioms together with $\neg F(\theta)$ for some term $\theta$ without variables. This gives an inconsistent theory, but one can study "concrete" consistency, in which proofs with only a "small" number of formulas cannot prove inconsistency.

To work within sequent calculus as before the preceding rules should be reformulated, but the point is clear enough and we omit the details. See also [12], [18].

Given a particular number $n$ one can prove $F(n)$ using the $F : successor$ rule $n$ times. Of course one can be more clever than that, using addition and multiplication to reduce substantially the length of the proof of the feasibility of $F(n)$. It is easy to imagine how one can get down to the realm of $\log n$, but on average one should not expect to do much better than that. We want to consider some particular large numbers for which one can find much shorter proofs of feasibility.

Define $e_2(n, 2)$ recursively for nonnegative integers $n$ by $e_2(0, 2) = 1$, $e_2(n + 1, 2) = 2^{e_2(n,2)}$. This function has *nonelementary* growth, i.e., it grows faster than any finite tower of exponentials. We want to describe a proof of $F(e_2(n, 2))$ with $O(n)$ lines due to Solovay.

Fix an $n$ for which we want to prove $F(e_2(n, 2))$. Define auxiliary formulas $F_i$, $i = 0, 1, 2, \ldots$ as follows. We take $F_0(x) = F(x)$, and we define $F_i(x)$

recursively by
$$F_i(x) \quad \text{is} \quad (\forall z)(F_{i-1}(z) \supset F_{i-1}(z^x)).$$

Conceptually we can imagine each formula $F_i(x)$ as defining a finite set of nonnegative integers. For each $i$ these are the integers $x$ such that the set corresponding to $F_{i-1}$ is closed under the mapping $z \mapsto z^x$. We can think of these subsets as becoming "smaller" as $i$ increases, and even in an exponential way.

The main property of these formulas $F_i(x)$ is that
$$F_i(x), F_i(y) \to F_i(x * y)$$

is provable. Let us first explain this in words. This sequent says that if $F_i(x)$ and $F_i(y)$ both hold, then $F_i(x * y)$ holds too. The hypotheses imply that $(\forall z)(F_{i-1}(z) \supset F_{i-1}(z^x))$ and that $(\forall w)(F_{i-1}(w) \supset F_{i-1}(w^y))$ are true. Applying the second with $w = z^x$ we get that $F_{i-1}(z^x) \supset F_{i-1}((z^x)^y)$ is true for all $x$, and hence $F_{i-1}(z) \supset F_{i-1}((z^x)^y)$ is true because of the assumption on $x$. This is the same as $F_i(x * y)$, because $(z^x)^y = z^{x*y}$. It is not hard to formalize this as a proof in sequent calculus, but we leave the details as an exercise.

Using this property one can show that $\to F_i(2)$ is provable for each $i$.

Next we claim that

(1) $$F_0(2), F_1(2), \ldots, F_n(2) \to F_0(e_2(n+1, 2))$$

is provable. To see this we use the following building blocks:

(2) $$F_i(e_2(n-i, 2)), F_{i-1}(2) \to F_{i-1}(2_2^e(n-i, 2))$$

These sequents are easy to prove using the definitions. They can be combined through cuts to give (1).

We can now conclude that
$$\to F_0(e_2(n+1, 2))$$

is provable. This follows from (1) and the provability of $\to F_i(2)$ for all $i$, established before. It is not hard to see that our proof had a total of $O(n)$ lines. Of course this is the same as $\to F(e_2(n+1, 2))$.

The *nesting* of quantifiers in this proof is responsible for the strong compression. Recall that $F_i(x)$ is defined recursively by $(\forall z)(F_{i-1}(z) \supset F_{i-1}(z^x))$. Thus the number of quantifiers in $F_i(x)$ is 1 plus twice the number of quantifiers in $F_{i-1}(z)$ and therefore grows exponentially in $i$. More importantly, in $F_i(x)$ we have two occurrences of $F_{i-1}$ which lie within the scope of the same quantifier and which are linked through the implication. By recursion there is a similar structure within $F_{i-1}$ and so forth. This nesting of links is crucial for the compression and leads to rich dynamical structure. The links combine to form an exponential number of cycles within the logical flow graph of the proof. See [12].

There is a notion of *curvatures* associated to the nesting of links between formula occurrences. This notion is used in [13] to read a finitely presented group from the cycles coming from a proof. For the example above the distortion of a certain subgroup inside the group reflects the compression in the proof and the expansion associated to cut-elimination.

# 10 A geometric view

How can we really "see" what happens in the preceding examples? Consider the John-Nirenberg theorem. The Calderón-Zygmund lemma is applied repeatedly to its own output to provide the existence of certain families of intervals without listing them explicitly. The proof of the feasibility of $e_2(n+1, 2)$ describes the number implicitly without writing out all the multiplications. In each of these cases an actual computation modelled on the proofs would have to visit certain formulas over and over again. There is cycling inside the proof with cuts.

To make this precise we use the concept of a *logical flow graph* introduced by Buss [8]. A different but related graph was introduced earlier by Girard [24]. Actually we shall use a modification of Buss' definition, in which we restrict ourselves to *atomic* formulas occurring in the proof. *Atomic formulas* are formulas without logical connectives. In propositional logic this means the propositional variables, while in predicate logic it means relations with their terms. The restriction to atomic formulas seems to be more natural for geometric and dynamical interpretations of proofs. Atomic formulas are like particles within a proof, and the logical flow graph traces their motion. To see cycling it is important to work with atomic formulas, to be able to move freely up and down the proof.

Let $\Pi$ be a proof of some sequent $S$. The vertices of the logical flow graph are the *occurrences* of atomic formulas in $\Pi$. We connect occurrences of atomic formulas by an edge only when they are *variants* of each other. In propositional logic two occurrences are variants when they represent the same formula, and in predicate logic we allow the terms within to be different.

We attach edges by tracing the logical relationships in a proof. Let us start with an axiom $A, \Gamma \to \Delta, A$. We do not attach any edges to the atomic subformulas in $\Gamma$ and $\Delta$. For the pair of distinguished occurrences of $A$ we attach an edge between each pair of their corresponding atomic formulas. For instance, if $A$ is $p \vee (\neg p \wedge q)$, then we attach edges between the two outer $p$'s, the two inner $p$'s, and the two $q$'s.

Next we have to decide how to attach edges in accordance with the rules. For example, consider the binary rule $\vee : left$,

$$\frac{A, \Gamma_1 \to \Delta_1 \quad B, \Gamma_2 \to \Delta_2}{A \vee B, \Gamma_{1,2} \to \Delta_{1,2}}$$

Each atomic formula on the top has a counterpart in the bottom, and we attach an edge between the two occurrences. No other edges are attached. Note that each atomic subformula in the $A$ on top is connected by an edge to its counterpart in the $A \vee B$ in the bottom, and similarly for $B$.

The other logical rules are treated in practically the same manner as the $\vee : left$ rule. The cut and contraction rules are treated differently. For the cut rule

$$\frac{\Gamma_1 \to \Delta_1, A \quad A, \Gamma_2 \to \Delta_2}{\Gamma_{1,2} \to \Delta_{1,2}}$$

we attach edges between the atomic subformulas of the two copies of the cut formula $A$, and between the atomic subformulas of the $\Gamma$'s, $\Delta$'s occurring in the upper sequents and their counterparts in the lower sequent. No other edges are attached, and in particular there are no edges from the cut formulas to the lower sequent. For a contraction,

$$\frac{\Gamma \to \Delta, A, A}{\Gamma \to \Delta, A} \qquad \frac{A, A, \Gamma \to \Delta}{A, \Gamma \to \Delta}$$

each atomic subformula in each of the two occurrences of $A$ in the upper sequent is connected to its counterpart in the lower sequent. Again the atomic occurrences in $\Gamma, \Delta$ in the upper sequent are connected to their counterparts in the lower sequent, and no other edges are attached. This is the only situation in which we have vertices attached to more than two edges.

Thus the cut and contraction rules are very different geometrically from logical rules. We shall not pursue this here, but there is much more to be said about this.

This completes the description of the logical flow graph, except for one extra ingredient, an *orientation*. Let us define first the *sign* of a formula. This simply counts the number of negations involved. Thus in the formula $p \wedge (\neg q \vee \neg(\neg r))$, the $p$ and $r$ occur positively and the $q$ occurs negatively. The logical connectives $\wedge, \vee, \forall, \exists$ do not affect the sign, but $\supset$ can. An atomic formula occurs positively in $A \supset B$ if it occurs positively in $B$ or negatively in $A$, and otherwise it occurs negatively. This is because of the implicit negation hidden in ⊃. Similarly, in a sequent $\Gamma \to \Delta$, an occurrence of an atomic formula $C$ is said to be *positive* in the sequent if it occurs positively in $\Delta$ or negatively in $\Gamma$, and it is said to be *negative* otherwise.

The logical flow graph has the nice property that any two variants connected by an edge have the same sign, except for axioms and cut rules, where we connect formulas of opposite sign.

We can use the signs of atomic formulas to define a natural orientation for the logical flow graph in the following manner. For an axiom we have edges from negative occurrences to positive occurrences. For logical rules and contractions edges between negative occurrences are oriented upwards and edges between positive occurrences are oriented downwards. In the cut rule we do the same except for the edges between the cut formulas, where the orientation goes from positive occurrences to negative occurrences. This orientation has a natural interpretation in logic [11].

This defines the logical flow graph associated to a proof as a directed graph. For the authors this kind of geometric picture is extremely helpful in trying to understand the structure of proofs. It provides a way to trace the logical relations between different occurrences of a formula in a sequent, and a more global alternative to the usual induction arguments.

To understand the meaning of the logical flow graph it is helpful to begin with some simple observations. Suppose that we have a proof $\Pi$ in propositional logic in which the occurrences of a variable $p$ come in distinct connected components of the graph. Then we can rename all of the occurrences of $p$ in one component and still have a valid proof.

Sometimes occurrences in a sequent have to be linked within any proof of it. For example, in the tautology $p \vee \neg p \to q \vee \neg q$ the two $q$'s must have the same name and therefore must be connected inside the proof. The two $p$'s do not need to be linked.

In complicated proofs one should expect many links between pairs of occurrences in the endsequent. For example, if all cut-free proofs $\Pi$ of a sequent $S$ are of exponential size compared to the size of $S$, then there must be a pair of occurrences in $S$ which are connected in the logical flow graph by an exponential number of links. This is not hard to derive from the definitions. (One has to be careful about weak occurrences, but they can be treated as in [11].) Remember from Section 3 that the cut-free proofs of the finite versions of the pigeon-hole principle in propositional logic must have exponential size.

From this we see that in some cases proofs must have *cycles*, i.e., paths in the logical flow graph which come back to where they started. However the cycles obtained as above are *unoriented*.

A basic observation (see [11]) is that proofs without cuts have no *oriented*

cycles. The reason for this is that in a proof without cuts we cannot have edges from positive occurrences to negative occurrences, and an oriented cycle must contain such an edge.

Proofs without contractions cannot have oriented cycles either. See [11, 13]. This is more amusing geometrically, one has to go further into the logical structure of the proof.

It is the combination of cuts and contractions that can lead to cycles. The proof described in Section 9 has an exponential number of cycles, as we mentioned before. Another example is given in [12] where cycles help to make a proof short, although less dramatically.

In general, the compression of proofs appears to be closely related to the presence of cyclic structures. This is reminiscent of certain phenomena of distortion in finitely presented groups. It is natural to try to use groups to represent dynamics within proofs, and this is explored in [13].

There are many natural questions concerning the relationship between cycles and the lengths of proofs. For instance, in [13] the question is raised of whether one can eliminate cyclic structures more efficiently than through cut elimination.

If we are going to admit the idea of geometric structures attached to proofs (a prominent theme in the algebraic approach to structure of proofs of Girard [28, 26, 27, 29]), then we should also think about mappings between spaces. This opens up a story too long for the present paper, but let us mention a basic point. The concept of mappings in the context of proofs suggests that we look for a more general notion of subproof which corresponds to embeddings. *Inner proofs*, introduced in [11], provide a candidate for such a notion. Roughly speaking, an inner proof can spread throughout the original proof (unlike a subproof), but keeping only some of the wires inside.

Inner proofs provide a notion of localization in a proof. We can ask about the local and global aspects of cut elimination. Given a proof $\Pi$ and a cut-free version $\Pi'$, one would like to know whether inner proofs in $\Pi'$ correspond to inner proofs in $\Pi$. This would be useful for matters of complexity, because it is easier to find inner proofs inside cut-free proofs, but then one would like to compress them as much as $\Pi$ is compressed. The usual method of cut-elimination does not work well for this endeavor, and an alternate method is introduced in [11] for propositional logic, where some conditions for making this transformation are given.

A basic point is that the known procedures of cut elimination either deform paths in the logical flow graph without breaking them, or they split the paths. One cannot recognize which alternative occurs in a local way. The problem of paths being split under cut elimination creates an obstruction for being able to transfer inner proofs in general. The method of [11] exploits the extra "determinism" in the history of propositional formulas to preserve more of the structure of the proof. In this method the splitting of paths is delayed until the elimination of atomic cuts. This can be accomplished by an explicit control of contractions.

These issues are related to Craig interpolation in a nice way. Not through the usual formulation of finding an interpolating formula, but in terms of splitting the given proof $\Pi : A \to B$ when given a truth assignment $\sigma$ to the common variables of $A$ and $B$ as mentioned at the end of Section 6. If one could always transfer inner proofs from the cut-free proof $\Pi'$ back to $\Pi$, then one could get linear bounds for the split proofs for a given $\sigma$. This would still not be enough to give an interpolating function (in the sense of Section 6) computable in polynomial time. All this remains open.

At this stage one starts to ask oneself a lot of questions about structure and proofs. Everything goes on the table, nothing is forbidden.

*IHES*
*35 Route de Chartres*
*91440 Bures-sur-Yvette*
*France*

A. Carbone
*Institut für Algebra und Diskrete Mathematik*
*Technische Universität Wien*
*Wiedner Hauptstrasse 8-10/118*
*A-1040 Wien*
*Austria*

S. Semmes
*Department of Mathematics*
*Rice University*
*Houston Texas 77251*
*U.S.A.*